\documentclass[conference]{IEEEtran}
\usepackage{amsmath,amssymb,amsthm,epsfig,color,subfigure,empheq,graphicx,graphics}
\usepackage{enumerate,url,algpseudocode,algorithm,wasysym,epstopdf,balance,enumitem}
\usepackage{multirow} 

\DeclareMathOperator{\nullspace}{null}

\newcommand \bzero{\mathbf{0}}
\newcommand \bone{\mathbf{1}}

\newcommand \be{\mathbf{e}}
\newcommand \bg{\mathbf{g}}

\newcommand \bp{\mathbf{p}}
\newcommand \bq{\mathbf{q}}

\newcommand \bu{\mathbf{u}}
\newcommand \bv{\mathbf{v}}
\newcommand \bw{\mathbf{w}}
\newcommand \bx{\mathbf{x}}

\newcommand \bA{\mathbf{A}}
\newcommand \bB{\mathbf{B}}
\newcommand \bC{\mathbf{C}}
\newcommand \bD{\mathbf{D}}

\newcommand \bF{\mathbf{F}}
\newcommand \bG{\mathbf{G}}

\newcommand \bI{\mathbf{I}}
\newcommand \bJ{\mathbf{J}}

\newcommand \bR{\mathbf{R}}

\newcommand \bX{\mathbf{X}}

\newcommand \btheta{\boldsymbol{\theta}}

\newcommand \blambda{\boldsymbol{\lambda}}

\newcommand \tbe{\tilde{\mathbf{e}}}

\newcommand \tbC{\tilde{\mathbf{C}}}
\newcommand \tbD{\tilde{\mathbf{D}}}

\newcommand \tblambda{\tilde{\boldsymbol{\lambda}}}

\newcommand \tbLambda{\tilde{\mathbf{\Lambda}}}

\newcommand \hbx{\hat{\mathbf{x}}}

\newcommand \hbJ{\hat{\mathbf{J}}}

\newcommand \bbe{\bar{\mathbf{e}}}

\newcommand \bbq{\bar{\mathbf{q}}}

\newcommand \bbC{\bar{\mathbf{C}}}
\newcommand \bbD{\bar{\mathbf{D}}}

\newcommand \bblambda{\bar{\boldsymbol{\lambda}}}

\def\BibTeX{{\rm B\kern-.05em{\sc i\kern-.025em b}\kern-.08em
		T\kern-.1667em\lower.7ex\hbox{E}\kern-.125emX}}

\begin{document}
	
\title{Learning to Optimize Power Distribution Grids using Sensitivity-Informed Deep Neural Networks}
	
\author{
\IEEEauthorblockN{
Manish K. Singh, 
Sarthak Gupta, 
Vassilis Kekatos}
\IEEEauthorblockA{\textit{Bradley Dept. of ECE, Virginia Tech}\\Blacksburg, VA, 24061\\Emails: \{manishks,gsarthak,kekatos\}@vt.edu}
\and
\IEEEauthorblockN{Guido Cavraro and Andrey Bernstein}
\IEEEauthorblockA{\textit{National Renewanable Energy Laboratories}\\Golden, CO, 80401\\Emails: \{guido.cavraro,Andrey.Bernstein\}@nrel.gov}
}

\maketitle
	
\begin{abstract} 
Deep learning for distribution grid optimization can be advocated as a promising solution for near-optimal yet timely inverter dispatch. The principle is to train a deep neural network (DNN) to predict the solutions of an optimal power flow (OPF), thus shifting the computational effort from real-time to offline. Nonetheless, before training this DNN, one has to solve a large number of OPFs to create a labeled dataset. Granted the latter step can still be prohibitive in time-critical applications, this work puts forth an original technique for improving the prediction accuracy of DNNs by taking into account the sensitivities of the OPF minimizers with respect to the OPF parameters. By expanding on multiparametric programming, it is shown that although inverter control problems may exhibit dual degeneracy, the required sensitivities do exist in general and can be computed readily using the output of any standard quadratic program (QP) solver. Numerical tests showcase that sensitivity-informed deep learning can enhance prediction accuracy in terms of mean square error (MSE) by 2-3 orders of magnitude at minimal computational overhead. Improvements are more significant in the small-data regime, where a DNN has to learn to optimize using a few examples. Beyond multiparametric QPs, the approach is currently being generalized to parametric (non)-convex optimization problems.
\end{abstract}
	
	
\section{Introduction}
\allowdisplaybreaks
\renewcommand{\thefootnote}{\fnsymbol{footnote}}
\footnotetext{Work partially supported by the US NSF grant 1751085.} 
To cope with resource allocation problems with stringent latency requirements, recent research advocates using learning models that are trained to solve optimization tasks; see e.g.,~\cite{eisen19}, \cite{ChenZhang20}. The general workflow involves solving a large number of such problems offline to create a dataset of problem instances paired with their minimizers; train a learning model to predict those minimizers; and use the model for promptly finding near-optimal solutions in real time. References \cite{eisen19} and \cite{Fioretto1} impose feasibility using Lagrangian decomposition, while \cite{Sidiropoulos18} and \cite{ZhangWangGiannakis19} build DNNs with structures mimicking algorithmic iterates solving the original optimization.  

In the context of power systems, there has recently been an interesting interplay between deep learning and the OPF. DNNs have been employed to predict OPF solutions under a linearized DC~\cite{DeepOPFPan19}, \cite{ZamzamBaker19}; or the exact AC model~\cite{Fioretto1}, \cite{GuhaACOPF}. The standard regression loss function is sometimes augmented by a penalty on OPF constraint violations to promote feasibility~\cite{DeepOPFPan19}, \cite{GuhaACOPF}, \cite{Fioretto1}. Alternatively, predicted solutions may be converted to a physically implementable schedule upon projection using a power flow (PF) solver~\cite{DeepOPFPan19}, \cite{ZamzamBaker19}. A graph neural network leveraging the connectivity of the power system is trained to infer AC-OPF solutions in~\cite{RibeiroGNNOPF}. Focusing on OPF for inverter control, reference~\cite{BaosenICNNvoltage} constructs a DNN to be used as a digital twin of the electric feeder. To learn optimal inverter control policies, DNN-based reinforcement learning schemes have been suggested in~\cite{ZhangISU},~\cite{Nanpeng19}.

Instead of learning minimizers, DNNs have also been used to predict the active constraints of an OPF~\cite{GuhaACOPF}, \cite{DekaMisraPowerTech19}, \cite{ChenZhang20}. References \cite{GuhaACOPF} and \cite{DekaMisraPowerTech19} train DNN-based classifiers to predict individual or sets of active constraints. In~\cite{ChenZhang20}, the active constraints of a linear programming (LP) DC-OPF are unveiled in three steps: A DNN is first trained to predict the optimal generation cost given demands; locational prices are then computed as the DNN sensitivities with respect to demands; and congested lines are identified via sparse-aware regression.

When the OPF is posed using a linearized grid model, it gives rise to parametric LPs or QPs, and hence, powerful tools from multiparametric programming (MPP) become relevant~\cite{BBM03}. Such tools have been utilized before to study congestion patterns in electricity markets~\cite{ZTL11}; for the probabilistic characterization of primal/dual OPF solutions~\cite{JTT17}; to accelerate hosting capacity analyses~\cite{TJKT20}; or strategic investment~\cite{TKV20}. MPP-aided methods however are meaningful only when the related OPF enjoys few congestion patterns and/or the application of interest is of a batch nature.

This work contributes to the literature of training DNNs to solve the OPF in the following creative way. A DNN is trained to fit not only the OPF minimizer, but also its sensitivities with respect to the problem parameters. Inspired by~\cite{Raissi19}, we cross-pollinate the idea of training sensitivity-aware DNNs from solving differential equations to learning for optimization problems in Section~\ref{sec:problem}. Advances in automatic differentiation facilitate the efficient computation of DNN sensitivities as well as their integration into training. In Sections~\ref{sec:dispatch} and \ref{sec:mpp}, we particularize properties from MPP theory to the inverter dispatch task at hand, so that the sensitivities of the OPF solution can be computed readily. The numerical tests of Section~\ref{sec:tests}
showcase that the proposed sensitivity-informed DNN outperforms a plain DNN in terms of prediction accuracy by 2-3 orders of magnitude. The paper is concluded with generalizations of the method to more challenging OPF tasks.

Regarding \emph{notation}, lower- (upper-) case boldface letters denote column vectors (matrices). Sets are represented by calligraphic symbols. Symbol $^\top$ stands for transposition, and inequalities are understood element-wise. All-zero and all-one vectors and matrices are represented by $\bzero$ and $\bone$; the respective dimensions are deducible from context. A symmetric positive definite matrix is denoted as $\bX\succ \mathbf{0}$.

\section{Sensitivity-Informed Deep Learning}\label{sec:problem}
Suppose a system operator has to routinely solve a convex quadratic program (QP) over the resource vector $\bx$
\begin{subequations}\label{eq:opf2}
	\begin{align}
	\bx_{\btheta}:=\arg\min_{\bx}\ &~\tfrac{1}{2}\bx^\top\bA\bx - \bx^\top\bB\btheta\label{eq:opf2:cost}\\
	\mathrm{s.to}\ &~\bC\bx\leq \bD\btheta+ \be: \quad \blambda_{\btheta}.\label{eq:opf2:c}
	\end{align} 
\end{subequations}
While $(\bA\succ \bzero,\bB,\bC,\bD,\be)$ are kept fixed for some time, the parameter vector $\btheta$ may be varying frequently. Let $(\bx_{\btheta},\blambda_{\btheta})$ denote a pair of primal/dual solutions for a particular $\btheta$. This multiparametric QP (MPQP) may be encountered in several resource allocation tasks in cyber-physical systems in general, and power systems in particular. Renditions of \eqref{eq:opf2} may appear in real-time electricity markets and contingency analysis in transmission systems; or in load disaggregation, demand response, and DER management in distribution systems. 

The application of interest here is dispatching inverters through an approximate OPF that minimizes power losses subject to voltage constraints. Vector $\bx$ comprises the inverter reactive power setpoints, while $\btheta$ carries the loading conditions (demand and solar generation) across feeder buses. Since $\btheta$ may change rapidly (rooftop solar generation can fluctuate by up to 15\% of its rating within one minute~\cite{Barker2012sustkdd}), an operator may not be able to solve \eqref{eq:opf2} exactly for thousands of feeders hosting thousands of inverters each. 

To expedite the process of computing inverter setpoints, the operator may want to train a DNN that \emph{learns to solve \eqref{eq:opf2}}. This DNN would be fed with $\btheta$ to compute a predictor $\hbx(\btheta;\bw)$ of $\bx_{\btheta}$. The straightforward approach for training this DNN is finding its weights $\bw$ as the minimizers of
\begin{equation}\label{eq:p-dnn}
\min_{\bw}\sum_{s=1}^S\|\hat{\bx}(\btheta_s;\bw)-\bx_{\btheta_s}\|_2^2.
\end{equation}
Problem \eqref{eq:p-dnn} is solved offline and aims at fitting the DNN output to the OPF minimizer over $S$ loading scenarios $\{\btheta_s\}_{s=1}^S$, which are deemed representative of the upcoming control period. Before solving \eqref{eq:p-dnn} however, the operator has to solve $S$ instances of \eqref{eq:opf2} to build the labeled dataset $\{(\btheta_s,\bx_{\btheta_s})\}_{s=1}^S$. Accommodating a large $S$ could be impractical in time-critical setups, where the DNN has to be re-trained every 30--60~min or so because parameters $(\bA,\bB,\bC,\bD,\be)$ change too. This is relevant for inverter control due to grid topology reconfigurations, changes in regulator settings and capacitor statuses, and varying loading conditions that may alter the linearization point of the nonlinear power flow equations.


We improve on this process in two aspects: \emph{i)} we reduce $S$ so the operator has to solve fewer OPFs; and \emph{ii)} for each OPF instance solved during training, we further exploit the sensitivity of $\bx_s$ with respect to $\btheta_s$. Both aspects can expedite training or increase the training dataset over the same training time, while \emph{ii)} is expected to yield better generalization to the sought DNN and stabilize its solutions against variations in $\btheta$. Inspired by \cite{Raissi19} where a DNN was trained to solve differential equations by incorporating information on its derivatives, here we integrate sensitivity information of the optimizer $\bx_{\btheta}$ with respect to the DNN input $\btheta$. To accomplish that, we propose augmenting the training process of \eqref{eq:p-dnn} as 
\begin{equation}\label{eq:si-dnn}
\min_{\bw}\sum_{s=1}^S\|\hat{\bx}(\btheta_s;\bw)-\bx_{\btheta_s}\|_2^2 + \|\hbJ(\btheta_s;\bw)-\bJ_{\btheta_s}\|_F^2
\end{equation}
where $\bJ_{\btheta}:=\nabla_{\btheta}\bx_{\btheta}$ and $\hbJ(\btheta;\bw):=\nabla_{\btheta}\hat{\bx}(\btheta;\bw)$ are the Jacobian matrices of the minimizer and the DNN output accordingly with respect to the parameter vector $\btheta$. Symbol $\|\cdot\|_F$ denotes the Frobenius norm. The optimization of \eqref{eq:si-dnn} aims at fitting the values of the minimizers as well as their sensitivities (partial derivatives) in terms of the DNN input $\btheta$. By incorporating such physics-aware information, the DNN is expected to predict well not only on the training scenarios $\btheta_s$'s, but in a neighborhood around them too. We term the DNN trained by \eqref{eq:si-dnn} as \emph{sensitivity-informed DNN} (SI-DNN) versus the \emph{plain DNN} (P-DNN) obtained from \eqref{eq:p-dnn}.


Having posed the sensitivity-informed deep learning task of \eqref{eq:si-dnn}, the pertinent questions are if the OPF operator $\mathrm{OPF}:\btheta\rightarrow \bx_{\btheta}$ is differentiable; and if so, whether its Jacobian matrix $\bJ_{\btheta}$ can be found in a computationally efficient way. The goal is to be able to build the augmented training set $\{(\btheta_s,\bx_{\btheta_s},\bJ_{\btheta_s})\}$ at a minimal computational effort in excess to the effort for building the training set $\{(\btheta_s,\bx_{\btheta_s}\}$ for \eqref{eq:p-dnn}. To this end, Section~\ref{sec:dispatch} first poses optimal inverter dispatching in the form of the MPQP in \eqref{eq:opf2}, and then Section~\ref{sec:mpp} adapts results from multiparametric programming to compute $\bJ_{\btheta}$ efficiently.

\section{Optimal Inverter Dispatch}\label{sec:dispatch}
Consider the task of optimal reactive power compensation by smart inverters. Given (re)active loads along with the active solar generation on every bus, the goal is to find the inverter setpoints for reactive power injections. Such setpoints can be found through an OPF that minimizes ohmic power losses on distribution lines, while satisfying voltage regulation limitations across buses. Before formally posing this OPF, we briefly review the postulated feeder model.

In a single-phase radial feeder with $N+1$ buses, the substation is indexed by $n=0$. Let $v_n$ be the voltage magnitude, and $p_n+jq_n$ the complex power injection at bus $n$. Collect all but the substation injections and voltages in the $N$-length vectors $(\bp,\bq,\bv)$. Vectors $(\bp,\bq)$ can be decomposed into non-negative injections $(\bp_i,\bq_i)$ from inverters, and withdrawals $(\bp_\ell,\bq_\ell)$ from loads as
\begin{subequations}\label{eq:pq}
\begin{align}
\bp&=\bF_i\bp_i-\bF_\ell\bp_\ell\label{eq:pq:p}\\
\bq&=\bF_i\bq_i-\bF_\ell\bq_\ell.\label{eq:pq:q}
\end{align}
\end{subequations} 
Matrices $\bF_i$ and $\bF_\ell$ map buses with solar and loads to the $N$ feeder buses, respectively. If there are $I$ buses hosting solar and $L$ buses hosting loads, then $\bF_i\in\{0,1\}^{N\times I}$ and $\bF_\ell\in\{0,1\}^{N\times L}$. We assume that each bus may host at most one (aggregate) load and at most one solar generator.

To simplify the task of dispatching inverters, we resort to an approximate feeder model obtained upon linearizing the power flow equations at the flat voltage profile of unit magnitudes and zero angles at all buses; see e.g., \cite{TJKT20}. According to this model, voltages depend linearly on power injections as 
\begin{equation}\label{eq:voltage}
\bv(\bp,\bq) \simeq \bR\bp + \bX\bq + v_0\bone_N
\end{equation}
where the $N\times N$ matrices $(\bR,\bX)$ are symmetric positive definite, and depend on the feeder topology and line impedances. We are also interested in expressing ohmic losses with respect to nodal injections. Adopting a second-order Taylor's series expansion, active losses on lines can be approximated as a convex quadratic function of power injections as
\begin{equation}\label{eq:losses}
L(\bp,\bq)\simeq 2\bp^\top\bR\bp + 2\bq^\top\bR\bq.
\end{equation}

Based on the previous modeling, the reactive setpoints for inverters can be found as the minimizers of
\begin{subequations}\label{eq:opf}
	\begin{align}
	\min_{\bq_i}\ &~\bq^\top\bR\bq\label{eq:opf:cost}\\
	\mathrm{s.to}\ &~-0.03\cdot\bone_N\leq \bR\bp+\bX\bq\leq 0.03\cdot\bone_N\label{eq:opf:c1}\\
	           &~-\bbq_i\leq \bq_i\leq \bbq_i.\label{eq:opf:c2}
	\end{align} 
\end{subequations}
The approximate OPF of~\eqref{eq:opf} minimizes ohmic losses over the controllable inverter setpoints $\bq_i$; note the first summand of \eqref{eq:losses} and the scaling by $2$ have been ignored. Constraint~\eqref{eq:opf:c1} maintains voltages within $\pm3\%$~per unit (pu). Although the substation voltage $v_0$ has been set at $1$~pu, it can be appended as an optimization variable of \eqref{eq:opf} without loss of generality. Our methodology can incorporate voltage regulators as in~\cite{TJKT20}; they are ignored here to keep the exposition succinct. Inverter setpoints are constrained in \eqref{eq:opf:c2} by the vector $\bbq_i$ of reactive ratings, which is assumed known and fixed.

Under varying grid conditions, an operator may want to solve~\eqref{eq:opf} routinely for different OPF instances. If the topology remains fixed, OPF instances differ in the values of loads and solar $(\bp_i,\bp_\ell,\bq_\ell)$. Define the vector of OPF inputs $\btheta:=[\bp^\top~\bq_\ell^\top]^\top$. Due to the problem structure, active loads and solar can appear collectively as $\bp$ per \eqref{eq:pq:p}. To better distinguish the optimization variable $\bx=\bq_i$ from the OPF parameters $\btheta$, let us substitute \eqref{eq:pq:q} in \eqref{eq:opf}, and express the inverter dispatch of \eqref{eq:opf} as the MPQP of \eqref{eq:opf2} with the substitutions
\begin{subequations}\label{eq:def}
	\begin{align}
	&\bA :=2\bF_i^\top \bR\bF_i \succ \bzero, \quad\quad
	\bB:=2\bF_i^\top\bR\bF_\ell\left[\bzero_{L\times N}~~\bI_L\right]\label{eq:def:AB}\\
	&\bC :=\left[\begin{array}{c}
	-\bX\bF_i\\
	\bX\bF_i\\
	-\bI_I\\
	\bI_I\end{array}\right],\quad\quad
	\bD:=\left[\begin{array}{c}
	\left[\bR~~-\bX\bF_\ell\right]\\
	\left[-\bR~~+\bX\bF_\ell\right]\\
	\bzero_{I}\\
	\bzero_{I}\end{array}\right]\label{eq:def:CD}\\
	&\be^\top :=\left[\begin{array}{c}
	0.03\cdot\bone_N^\top~~~
	0.03\cdot\bone_N^\top~~~
	\bbq_i^\top~~~
	\bbq_i^\top
	\end{array}\right].\label{eq:def:e}
	\end{align} 
\end{subequations}

\section{Computing the Jacobian Matrix via MPP}\label{sec:mpp}
Consider solving \eqref{eq:opf2} for a particular $\btheta$. The corresponding optimal primal-dual pair $(\bx_{\btheta},\blambda_{\btheta})$ should satisfy the Karush-Kuhn-Tucker (KKT) optimality conditions~\cite{Be99}. Among the inequality constraints of \eqref{eq:opf2}, those satisfied with equality at the optimum are termed \emph{active or binding constraints}. Let $(\tbC,\tbD,\tbe,\tblambda_{\btheta})$ be the row partitions of $(\bC,\bD,\be,\blambda_{\btheta})$ associated with active constraints. Likewise, let $(\bbC,\bbD,\bbe,\bblambda_{\btheta})$ denote the row partition related to \emph{inactive or non-binding constraints}. Leveraging complementarity slackness, the KKT conditions can be expressed as~\cite{Be99}
\begin{subequations}\label{eq:kkt}
	\begin{align}
	\bA\bx_{\btheta}+\tbC^\top\tblambda_{\btheta}&=\bB\btheta \label{eq:kkt:a}\\
	\tbC\bx_{\btheta}&=\tbD\btheta+\tbe \label{eq:kkt:b}\\
	\bbC\bx_{\btheta}&<\bbD\btheta+\bbe \label{eq:kkt:c}\\
	\tblambda_{\btheta}&> \bzero \label{eq:kkt:d}\\
	\bblambda_{\btheta}&= \bzero. \label{eq:kkt:e}
	\end{align} 
\end{subequations}
For \eqref{eq:kkt:d} in particular, it is assumed that binding constraints relate to strictly positive dual variables $\tblambda_{\btheta}>\bzero$. The reverse is unlikely to occur when $\btheta_s$'s are drawn at random. In the unlikely case some of the entries of $\tblambda_{\btheta}$ do equal zero, the OPF sensitivity for this specific $\btheta$ can be ignored and the DNN is trained only to match the OPF minimizer $\bx_{\btheta}$.

If there are $A$ active constraints, the equalities of \eqref{eq:kkt:a}--\eqref{eq:kkt:b} constitute a system of $(I+A)$ linear equations over the $(I+A)$ unknowns $(\bx_{\btheta},\tblambda_{\btheta})$. Since $\bA\succ \bzero$, the Lagrangian optimality condition of \eqref{eq:kkt:a} yields 
\begin{equation}\label{eq:x}
\bx_{\btheta}=\bA^{-1}(\bB\btheta-\tbC^\top\tblambda_{\btheta}).     
\end{equation}
Substituting \eqref{eq:x} into \eqref{eq:kkt:b} provides
\begin{equation}\label{eq:lambda}
\bG\tblambda_{\btheta}=(\tbC\bA^{-1}\bB-\tbD)\btheta -\tbe
\end{equation}
where $\bG:=\tbC\bA^{-1}\tbC^\top$. If matrix $\tbC$ has linearly independent rows, then $\bG\succ \bzero$, and \eqref{eq:lambda} has a unique solution. Otherwise, there are infinitely many $\tblambda_{\btheta}$'s satisfying \eqref{eq:lambda} within a shift invariance on the nullspace $\nullspace(\tbC^\top)=\nullspace(\bG)$.

\begin{figure}[t]
    \centering
    \hspace*{.5em}\includegraphics[scale=0.8]{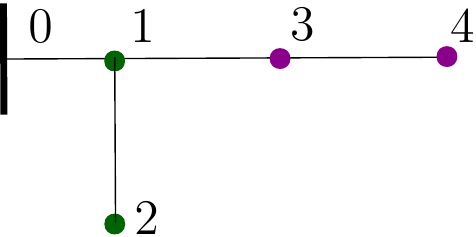}
    \caption{A 5-bus feeder demonstrating the violation of LICQ for \eqref{eq:opf}.}
    \label{fig:5bus}
\end{figure}

For a general QP, the rows of $\tbC$ are typically linearly independent, thus satisfying the so termed linear independence constraint qualification (LICQ)~\cite{Be99}. Nonetheless, that is not the case for the inverter dispatch problem of \eqref{eq:opf} where instances of linearly dependent active constraints occur frequently. To see this, consider the feeder of Fig.~\ref{fig:5bus}, where buses $\{1,2\}$ host inverters, bus $3$ is a load bus, and $4$ is a zero-injection bus. Let $x_n$ denote the reactance for the line feeding bus $n$. Then, matrix $\bX$ can be found as~\cite{TJKT20}
\begin{equation*}
\bX=\begin{bmatrix}
x_1&x_1&x_1&x_1\\
x_1&x_1+x_2&x_1&x_1\\
x_1&x_1&x_1+x_3&x_1+x_3\\
x_1&x_1&x_1+x_3&x_1+x_3+x_4
\end{bmatrix}.
\end{equation*}
Moreover, the product $\bX\bF_i$ selects the first two columns of $\bX$ corresponding to inverter locations. Now consider a loading scenario that causes the voltage at bus $3$ to be at the lower limit. Being a  zero-injection bus, the downstream bus $4$ shares the same voltage. Thus, the minimum voltage constraints become active for buses $3$ and $4$. Assuming no other constraints being active, matrix $\tbC$ of \eqref{eq:def:CD} reads as \[\tbC=-\begin{bmatrix}x_1&x_1\\x_1&x_1\end{bmatrix}.\]
which is clearly row rank-deficient. This counterexample emphasizes the need to explicitly address the case where the solution of \eqref{eq:opf2} or \eqref{eq:opf} fails to meet the LICQ.

Regardless of whether $\tbC$ is full row-rank or not, the solution to \eqref{eq:lambda} lies in the polyhedron
\begin{align*}
\tbLambda_{\btheta}:=\{\tblambda_{\btheta}=\bG^{+}(\tbC\bA^{-1}\bB-\tbD)\btheta -\bG^{+}\tbe +\bu,\tbC^\top\bu=\bzero\}
\end{align*}
where $\bG^{+}$ is the pseudo-inverse of $\bG=\tbC\bA^{-1}\tbC^\top$. Substituting $\tblambda_{\btheta}$ back in \eqref{eq:x} and exploiting $\tbC^\top\bu=\bzero$, we get
\begin{equation}\label{eq:x2}
\bx_{\btheta}=\bJ_{\btheta}\btheta +\bA^{-1}\tbC^\top\bG^{+}\tbe.    
\end{equation}
where the sensitivity matrix is defined as
\begin{equation}\label{eq:Jacobian}
\bJ_{\btheta}:=\bA^{-1}\bB-\bA^{-1}\tbC^\top\bG^{+}(\tbC\bA^{-1}\bB-\tbD).
\end{equation}
Although there may be infinitely many $\tblambda_{\btheta}$'s satisfying the KKT conditions, the optimal primal solution of \eqref{eq:opf2} is unique if it exists. This is not surprising since \eqref{eq:opf2} has a strictly convex objective. Moreover, the solution can be expressed as an affine function of $\btheta$. Note that the parameters of this affine function depend on the set of active constraints, and this is indicated by the subscript of $\btheta$ on $\bJ_{\btheta}$.

Because the triplet $(\btheta,\bx_{\btheta},\tblambda_{\btheta})$ should also satisfy \eqref{eq:kkt:d} and \eqref{eq:kkt:c}, there exists a $\bu\in\nullspace(\tbC^\top)$ so that $\btheta$ satisfies
\begin{subequations}\label{eq:cr}
\begin{align}
    \bG^{+}(\tbC\bA^{-1}\bB-\bD)\btheta >\bG^{+}\tbe -\bu\label{eq:cr:a}\\
    (\bbC\bJ_{\btheta}-\bbD)\btheta<\bbe-\bbC\bA^{-1}\tbC^\top\bG^{+}\tbe\label{eq:cr:b}
\end{align}    
\end{subequations}

So far, we have characterized the solution to \eqref{eq:opf2} for a particular $\btheta$. Since we are interested in the sensitivity of the OPF minimizer with respect to $\btheta$, we ask the natural question whether \eqref{eq:x2} holds for all $\btheta'$'s within a vicinity of this specific $\btheta$. The answer to this question is on the affirmative. To see this, fix $\bu$ and perturb $\btheta$ to get $\btheta'$ so it still satisfies \eqref{eq:cr}. Using $\btheta'$ in lieu of $\btheta$, construct $\bx_{\btheta'}$ from \eqref{eq:x2}, and $\tblambda_{\btheta'}$ from \eqref{eq:lambda}. In doing so, we row-partition matrices assuming the same constraints are active as for $\btheta$. This means we still use matrix $\bJ_{\btheta}$. We also set $\bblambda_{\btheta'}=\bzero$. The constructed primal-dual pair $(\bx_{\btheta'},\blambda_{\btheta'})$ satisfies the KKT conditions of \eqref{eq:opf2} for $\btheta'$, and hence constitutes an optimal solution for this $\btheta'$. The aforesaid process can be repeated for any $\btheta'$ in the vicinity of the original $\btheta$ because \eqref{eq:cr} are strict inequalities. In other words, formula \eqref{eq:x2} is valid for all $\btheta'$'s around $\btheta$ and with matrix $\bJ_{\btheta}$ remaining unaltered. 

The latter reveals that the OPF operator from $\btheta$ to $\bx_{\btheta}$ is differentiable around $\btheta$. In addition, its Jacobian matrix is provided in closed form using \eqref{eq:Jacobian}. Calculating $\bJ_{\btheta}$ entails: \emph{i)} Knowing the set of active constraints; \emph{ii)} inverting matrix $\bA$ once; and \emph{iii)} inverting matrix $\bG$. Although \emph{iii)} is executed once per $\btheta$, the computation is lightweight since the number of active constraints $A$ should be smaller than $N$. In a nutshell, computing $\bJ_{\btheta}$ once \eqref{eq:opf2} has been solved for a particular $\btheta$, is much simpler than solving \eqref{eq:opf2} per se. Hence, computing sensitivities adds insignificant complexity in the process of constructing the labeled dataset $(\btheta,\bx_{\btheta},\bJ_{\btheta})$.

\begin{figure}
    \centering
    \hspace*{.5em}\includegraphics[scale=0.22]{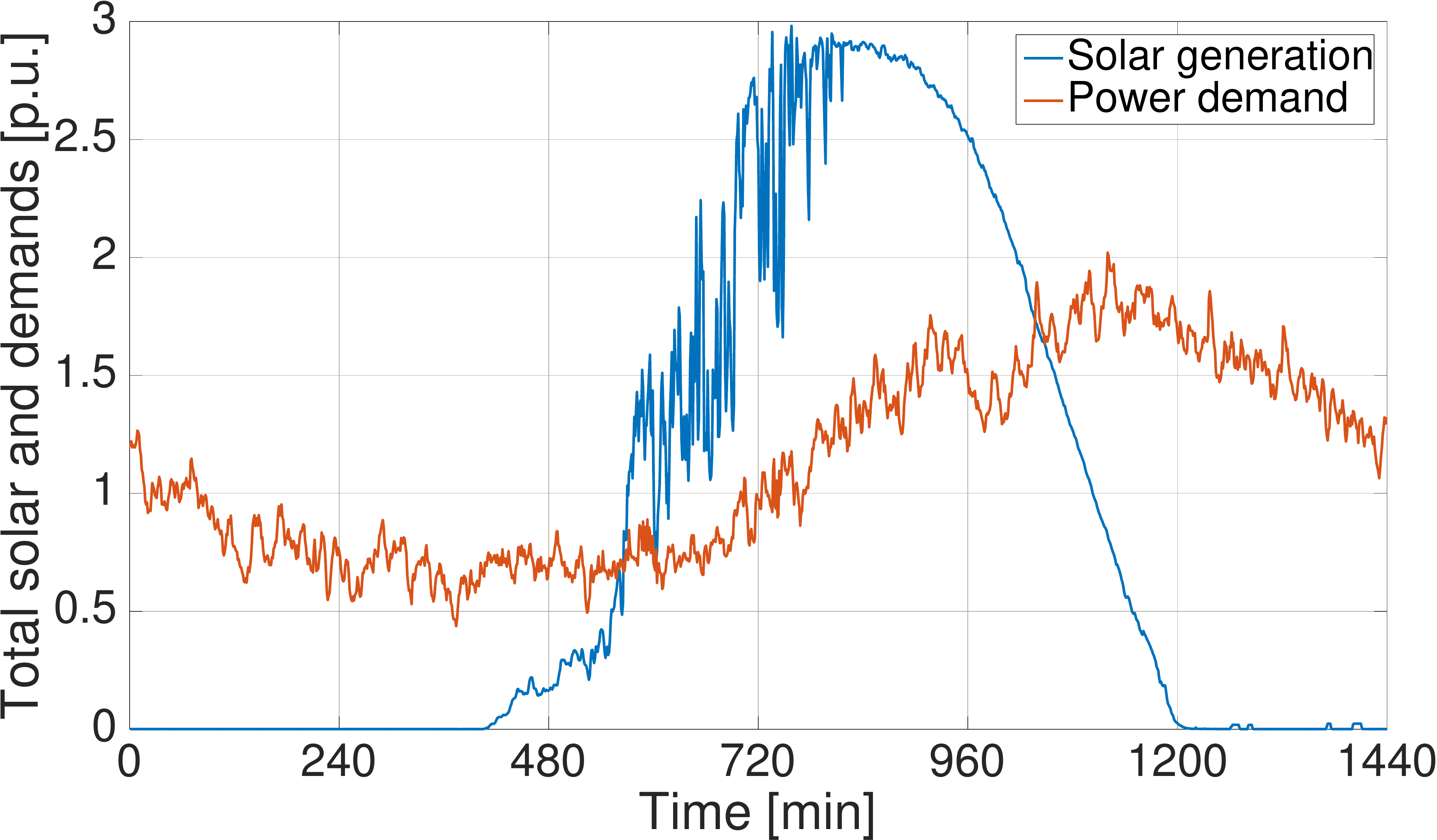}\\
    \includegraphics[scale=0.22]{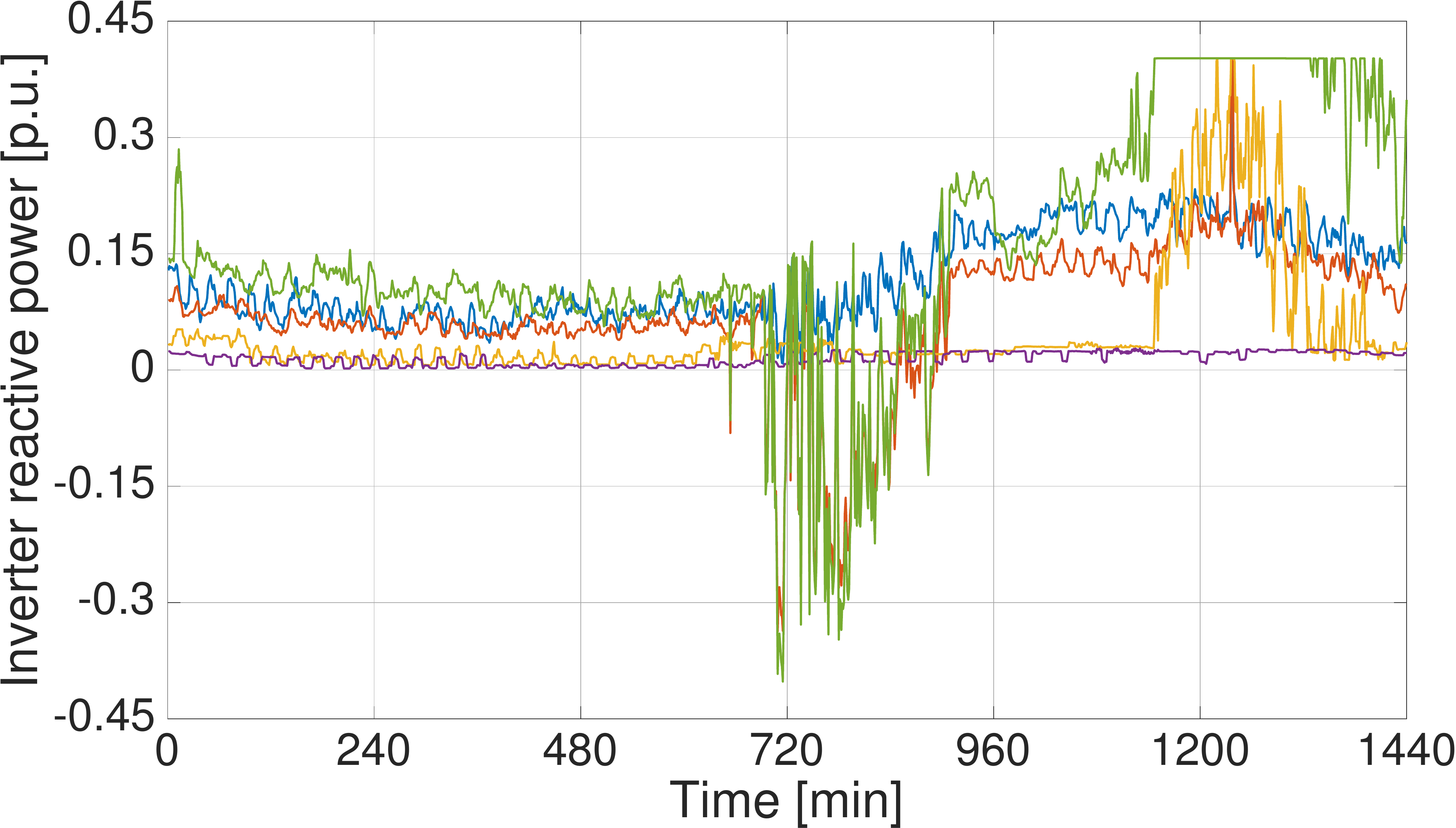}
    \caption{Minute-based data for the total (feeder-wise) solar power generation and active power demand (top panel); and the optimal reactive power injection setpoints for five inverters (bottom panel).}
    \label{fig:solarload}
\end{figure}

\section{Numerical Tests}\label{sec:tests}

\begin{figure}
    \centering
    \hspace*{.5em}\includegraphics[scale=0.29]{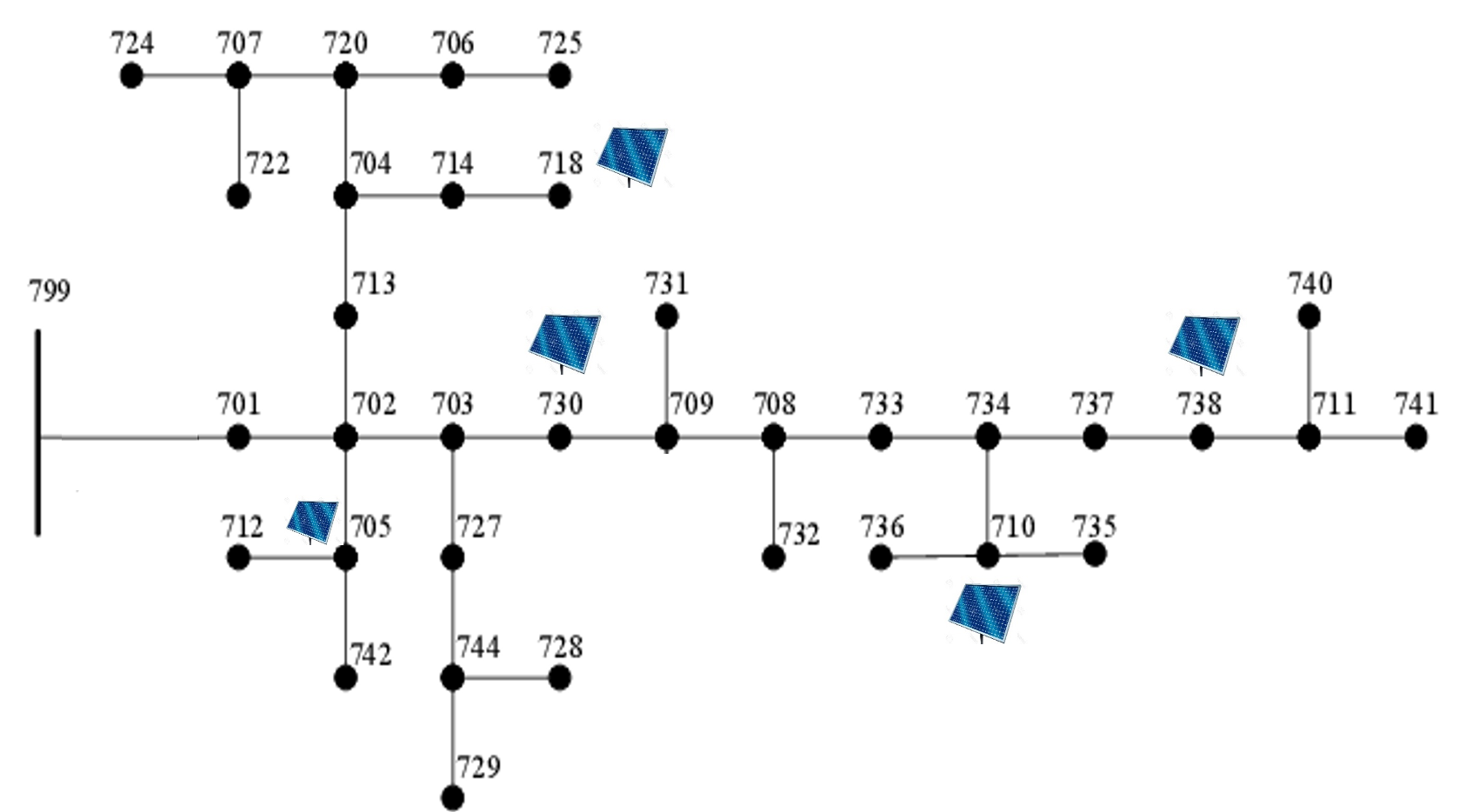}
    \caption{A modified IEEE 37-bus feeder showing additional solar generators.}
    \label{fig:37bus}
\end{figure}

The performance of the developed approach was numerically tested by solving the OPF in \eqref{eq:opf} for the IEEE 37-bus feeder upon removing regulators, incorporating 5 solar generators, and converting it to its single-phase equivalent; see Fig.~\ref{fig:37bus}. We extracted minute-based load and solar generation data for June 1, 2018, from the Pecan Street dataset~\cite{pecandata}. The feeder has 25 buses with non-zero load. The first 75 non-zero load buses from the dataset were aggregated every 3 and normalized to obtain 25 load profiles. Similarly, 5 solar generation profiles were obtained. The normalized load profiles for the 24-hr period were scaled so the 97-th percentile of the total load duration curve coincided with the total nominal load. This scaling results in a peak aggregate load being 1.1 times the total nominal load. We synthesized reactive loads by scaling active demand to match the power factors of the IEEE 37-bus feeder. The 5 solar generators were of equal ratings and scaled to meet 75\% of the total energy requirement over the entire day. Figure~\ref{fig:solarload} shows the total demand and solar generation across the feeder (top), along with the optimal reactive power injection for each of the five inverters (bottom). The optimal setpoints were obtained by solving~\eqref{eq:opf} using YALMIP and Sedumi. In solving the 1,440 OPF instances, no constraints were active for 914 instances. Out of the remaining 526 instances, the LICQ was not satisfied for 175 instances; thus necessitating the approach pursued here. 


Our tests compared the P-DNN and SI-DNN, both trained to predict the minimizer of \eqref{eq:opf}. For the first set of tests, the DNNs were assumed to be trained on an hourly basis. To evaluate the potential benefit of integrating sensitivity information, the architecture, optimizer, and training epochs were kept identical for the two DNNs. In fact, the architecture was chosen in favor of the P-DNN. It was determined so that for a sufficiently large training set, the P-DNN would be able to predict a near-optimal inverter dispatch. This ensured a sufficiently rich, yet not too complex architecture: Three hidden layers with 210, 210, and 350 neurons respectively, all using the rectified linear unit (ReLU) activation.  All results were generated using the Adam optimizer with a learning rate of~0.01. The DNNs were implemented using the TensorFlow library on Google colab.

\begin{figure}[t]
    \centering
    \includegraphics[scale=0.2]{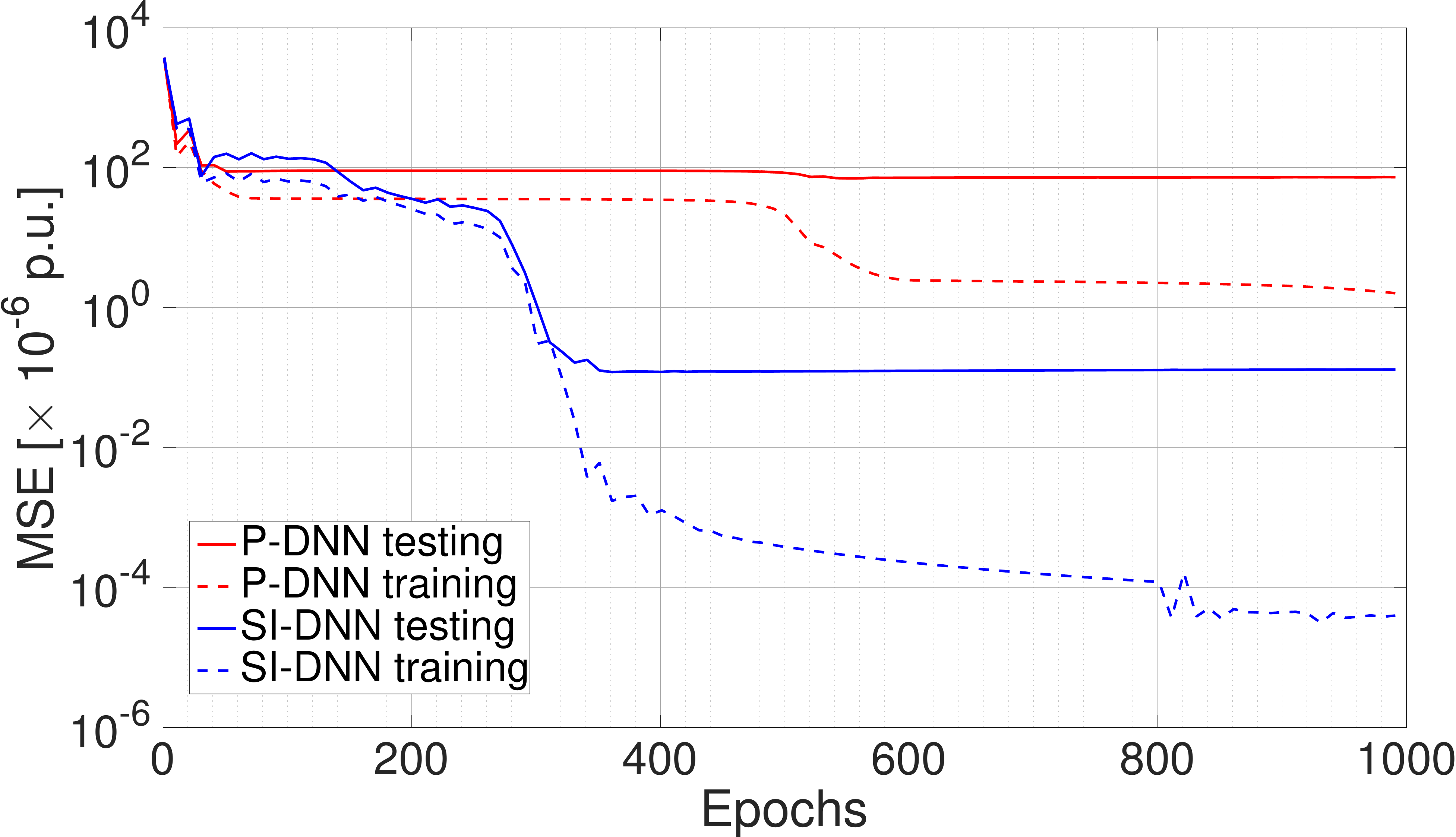}\vspace*{0.4em}
    \includegraphics[scale=0.2]{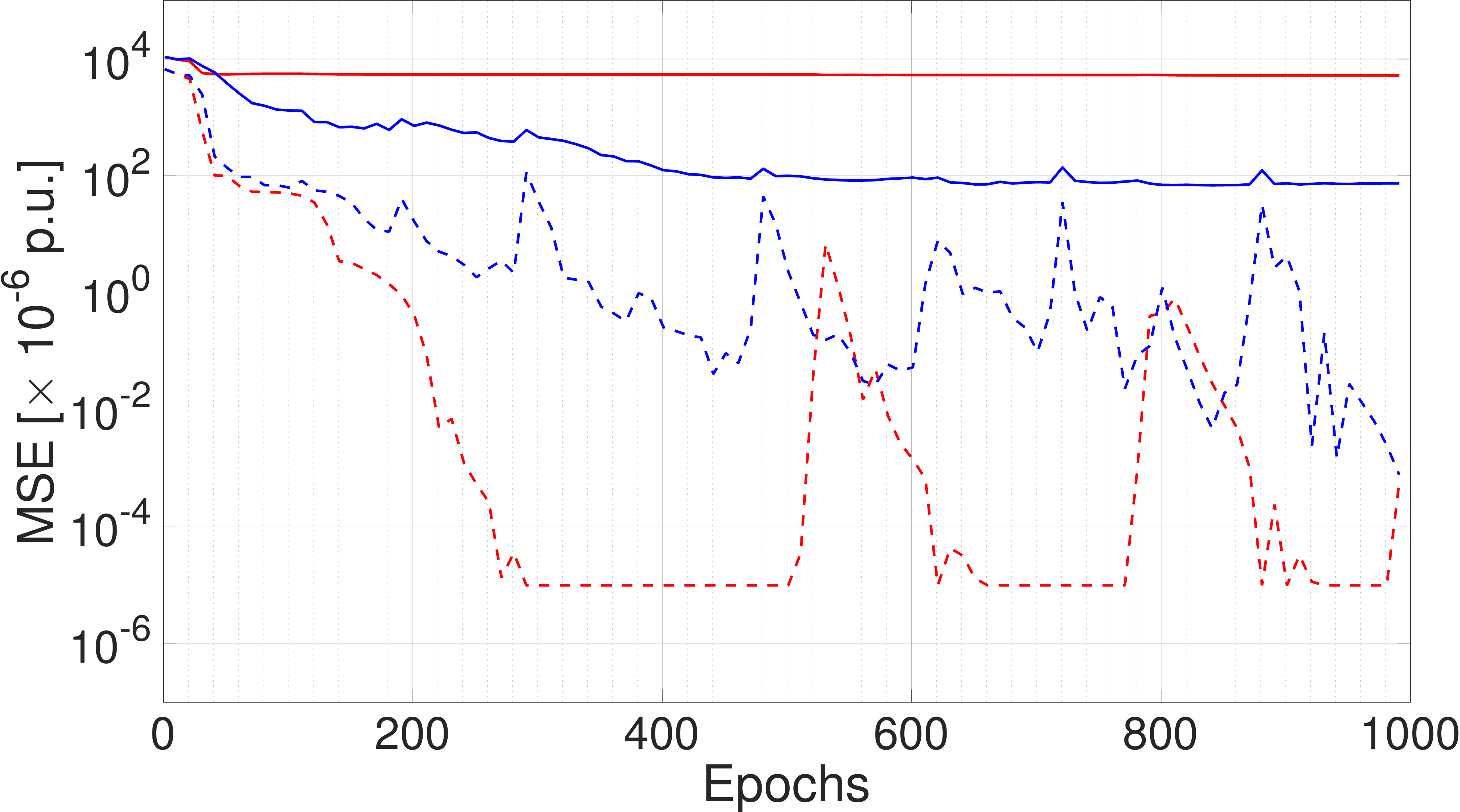}
    \includegraphics[scale=0.2]{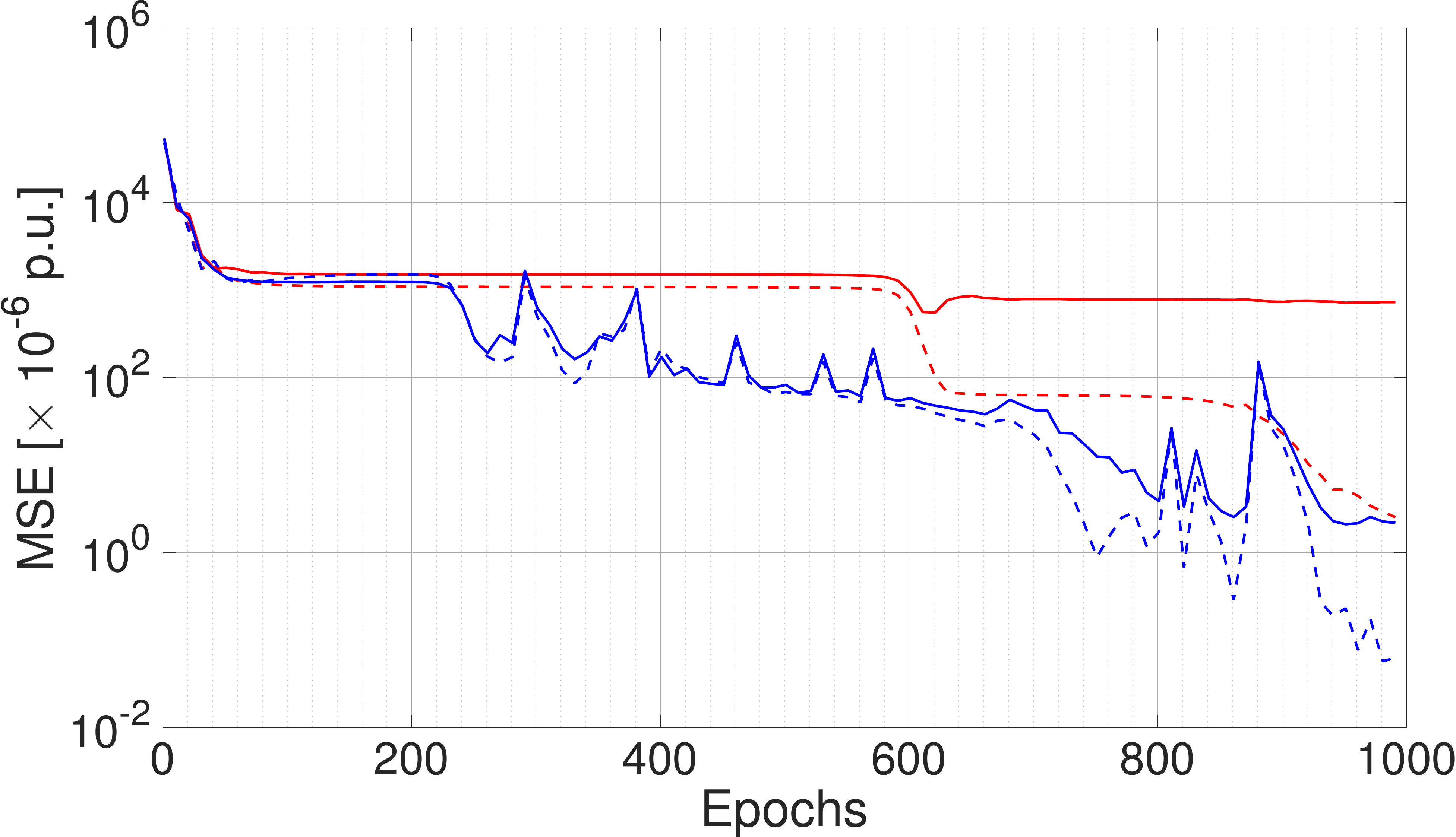}
    \caption{Training and testing errors across epochs in terms of MSE for hours 5 (top), 12 (middle), and 20 (bottom). Training was performed using 15-min smart meter data (4 OPF instances) to operate during the hour of interest.}
    \label{fig:hour5-12-20}
\end{figure}

For a given hour, the two DNNs were trained using 4 OPF instances, obtained by sampling the upcoming hourly period every 15~min. The remaining 56 grid instances were used to test the two DNNs based on the mean square error (MSE) $\|\bx_s-\hat{\bx}_s\|_2^2$ on the predicted setpoints. To compare P-DNN and SI-DNN under varying conditions, their training and testing errors were evaluated over hours 5, 12, and 20, and shown in Fig.~\ref{fig:hour5-12-20}. The tests show that even when P-DNN yields smaller training errors, SI-DNN offers an improvement in testing error by one or two orders of magnitude.

\begin{table}[t]
	\caption{Average Test MSE [in $10^{-6}$ pu] after 1,000 Epochs\newline for Different Hours of the Day}
	\begin{center}
		\begin{tabular}{c|rr|rr}
			\hline\hline
	        \multirow{2}{*}{Hour}&\multicolumn{2}{c|}{12 OPF training scenarios}&\multicolumn{2}{c}{4 OPF training scenarios}\\\cline{2-5}	        
	        & \textbf{P-DNN} & \textbf{SI-DNN} & \textbf{P-DNN} & \textbf{SI-DNN}\\
	        \hline\hline
	        5& 27.9&0.04&75.3&0.33\\
	        \hline
	        12& 1184.4 & 50.4 & 7229.8 & 4755.3 \\
	        \hline
	        16& 145.1 & 49.6 & 540.3 & 342.5 \\
	        \hline
	        20 & 359.3 & 45.0 & 1500.6 & 145.2 \\
			\hline\hline
			\end{tabular}
		\end{center}
		\label{tbl:testloss}
		\vspace*{-1.5em}
\end{table}

We next conducted the previous test but now over 10 Monte Carlo runs, each time selecting the training OPF scenarios at random. To evaluate the effect of the dataset size, we trained both DNNs using 4 and 12 scenarios for each of the hours 5, 12, 16, and 20. Table~\ref{tbl:testloss} reports the test MSEs attained after 1,000 epochs and averaged over the 10 Monte Carlo runs. The results corroborates that the SI-DNN features a gain in prediction accuracy over P-DNN by 1-3 orders of magnitude. Its advantage is generally more pronounced when dealing with smaller datasets, as expected.


\begin{table}[t]
	\caption{Average Test MSE [in $10^{-6}$ pu] and training Time [in sec]\newline after 1,000 Epochs for 10--12~am (840 Min-based Scenarios)}
	\vspace*{-1.5em}
	\begin{center}
		\begin{tabular}{c|rr|rr}
			\hline\hline
	        Training &\multicolumn{2}{c|}{\textbf{P-DNN}}&\multicolumn{2}{c}{\textbf{SI-DNN}}\\\cline{2-5}	        
	        Scenarios& \textbf{MSE} & \textbf{Time} & \textbf{MSE} & \textbf{Time}\\
	        \hline\hline
	        5\%& 1407.8 & 108.9 & 119.9 & 118.7 \\
	        \hline
	        10\%& 520.2& 77.1 & 72.9 & 79.7  \\
	        \hline
	        20\%& 219.2& 100.8 & 55.1 & 113.4  \\
			\hline\hline
			\end{tabular}
		\end{center}
		\label{tbl:testloss_fullday}
		\vspace*{-1.5em}
\end{table}

For the last set of tests, the two DNNs were trained to learn the OPF for an entire day. The training dataset was constructed by sampling 5\%, 10\%, and 20\% of the one-minute data. To explicitly focus on periods of high variability, we excluded the hours from midnight to 10~am from sampling. Table~\ref{tbl:testloss_fullday} summarizes the MSEs obtained after 1,000 epochs and averaged over 100 Monte Carlo draws of the training dataset. The table also reports the average training time for 1,000 epochs for the two methodologies. The runtimes on Google colab often varies with session. Thus, these times can only be compared separately for each training scenario; and not across scenarios. The results establish that the SI-DNN consistently outperforms the P-DNN without incurring any significant increase in training time. For example, the SI-DNN achieves an MSE of $119.9\cdot 10^{-6}$~pu using 5\% of the data, whereas the P-DNN achieves an MSE of $219.2\cdot 10^{-6}$~pu although it has been trained by using 4 times more data (20\%).



\section{Conclusions and Ongoing Work}\label{sec:conc}
A novel procedure for training DNNs that learn to optimize has been put forth. Leaping beyond the general practice of training DNNs via labeled parameter-minimizer pairs, this work ensures that the sensitivities of DNN predictions to inputs are close to the respective sensitivities of the original optimization task. Addressing QPs in particular, the sensitivities required for training the DNN can be computed readily in virtue of results from MPP theory. The application pursued has been inverter reactive power control in distribution systems for minimizing losses subject to voltage constraints. It has been shown here that although for general QPs dual degeneracies are rare, for the inverter dispatch task such degeneracies can be encountered frequently. Fortunately, even for such instances, the required Jacobian matrices do exist in general, and the proposed approach can successfully compute them in closed form. The efficacy of the novel training method is demonstrated via numerical tests that corroborate an improvement in prediction accuracy by 2-3 orders of magnitude, as compared to the traditional regression approach. The improvements are more pronounced in the small-data regime, where a DNN has to learn to optimize using few examples.

Prompted by these promising results, we are currently generalizing this creative idea of including gradient information into DNNs towards several exciting directions. Beyond MPQPs that feature rich structure in their solutions, we are interested in improving learning for resource allocation tasks of the generic parametric form
\begin{align}%
\bx_{\btheta}:=\arg\min_{\bx}~&~f(\bx;\btheta)\tag{$P_{\btheta}$}\\
\mathrm{s.to}~&~\bg(\bx;\btheta)\leq \bzero:~\blambda_{\btheta}.\nonumber
\end{align}
Thanks to results from optimization, given $(f,\bg,\bx_{\btheta},\blambda_{\btheta})$ one can compute the gradient $\nabla_{\btheta}f(\bx_{\btheta})$, and the Jacobian matrices $\nabla_{\btheta}\bx_{\btheta}$ and $\nabla_{\btheta}\blambda_{\btheta}$, regardless if $(P_{\btheta})$ is convex. Leveraging these results, we are currently working towards: \emph{d1)} incorporating sensitivities for learning to optimize different power system optimization and monitoring tasks, including the AC OPF and its various convex relaxations; \emph{d2)} integrating dual sensitivities and weight their deviations; and \emph{d3)} predicting binding constraints so that exact minimizers can be found by solving an OPF over a condensed feasible set.  
 
\balance
\bibliographystyle{IEEEtran}
\bibliography{myabrv,inverters}	
\end{document}